\documentclass{amsart}
\usepackage{amscd, amsmath, amsthm, amssymb, xypic, multirow}
\usepackage[all, knot]{xy}
\xyoption{arc}

\newtheorem{theorem}{Theorem}[section]
\newtheorem{lemma}[theorem]{Lemma}
\newtheorem{proposition}[theorem]{Proposition}

\theoremstyle{definition}     
\newtheorem{notation}[theorem]{Notation}

\theoremstyle{remark}
\newtheorem{remark}[theorem]{Remark}

\def \cen #1 { \begin{center} #1 \end{center}}

\begin{document}

\title[Connecting rigid Calabi--Yau threefolds]
{Connecting certain rigid birational non-homeomorphic Calabi--Yau
threefolds via Hilbert scheme}

\author[N.-H.\ Lee]{Nam-Hoon Lee}
\address{School of Mathematics, Korea Institute for Advanced Study,
Dongdaemun-gu, Seoul 130-722, Korea}
\email{nhlee@kias.re.kr}

\author[K. Oguiso]{Keiji Oguiso}
\address{Department of Economics, Keio University, Hiyoshi Kohoku-ku, Yokohama, Japan and School of Mathematics, Korea Institute for Advanced Study,
Dongdaemun-gu, Seoul 130-722, Korea
}
\email{oguiso@hc.cc.keio.ac.jp}
\thanks {The second named author was supported by JSPS}

\subjclass[2000]{14J32, 14D06}

\begin{abstract} We shall give an explicit pair of birational projective
Calabi--Yau threefolds which are rigid, non-homeomorphic, but are connected by
 projective flat deformation over some connected base scheme.
\end{abstract}

\maketitle


\setcounter{section}{-1}
\section{Introduction}
A \emph{Calabi--Yau manifold} is a compact K\"ahler
simply-connected manifold with a nowhere vanishing global $n$-form
but no global $i$-from with $0< i < n=\dim X $. By Kodaira's
criterion, it is projective if dimension $n \ge 3$.

As well known, Calabi--Yau manifolds, hyperk\"ahler manifolds and
complex tori form the building blocks of compact K\"ahler
manifolds with vanishing first Chern class (\cite{Be1},
\cite{GHJ}). A famous theorem of Huybrechts states that two
bimeromorphic hyperk\"ahler manifolds are equivalent under smooth
deformation (\cite{Hu1}, \cite{Hu2}). In particular, they are
homeomorphic to each other, having the same Betti numbers and
Hodge numbers. Clearly, the same holds true for complex tori.
Other famous theorem, called Kontsevich's theorem, says that two
birational Calabi--Yau manifolds have the same Betti numbers and
Hodge numbers (\cite{Ba}, \cite{DL}, \cite{It}, \cite{Ya},
\cite{Wa}). However, there are rigid birational non-isomorphic
Calabi--Yau manifolds (cf.\ Theorem \ref{mainthm}). Obviously they
are not equivalent under any \emph{smooth} deformation.

The aim of this note is to remark that there nevertheless exist
birational Calabi--Yau threefolds which are rigid,
non-homeomorphic,
but are connected by (necessarily non-smooth)
projective flat deformation:

\begin{theorem}\label{mainthm} There are  Calabi--Yau threefolds $X$ and $Y$ such that:
\begin{enumerate}
\item $X$ and $Y$ are birational and rigid,
\item $X$ and $Y$ are not homeomorphic but, \label{non-homeo}
\item $X$ and $Y$ are connected by projective flat deformation over some connected scheme.
\label{conn}
\end{enumerate}
\end{theorem}

This work is motivated by the famous fantasy of Miles Reid
\cite{Re} -- especially the question what this fantasy would be
like for rigid Calabi--Yau threefolds -- and by the first named
author's
recent result on the equivalence of certain Calabi--Yau
threefolds with Picard number one, of different topological type,
under projective flat deformation
\cite{Le}.

In the proof of our main theorem (Theorem \ref{mainthm}), the
following  deep theorem of Hartshorne \cite{Ha} (see also
\cite{PS}) plays an important role:
\begin{theorem}[R. Hartshorne]
A Hilbert scheme ${\rm Hilb}^{P(x)}_{\mathbf P^N}$ of $\mathbf P^N$
with fixed Hilbert polynomial $P(x)$ is connected.
\end{theorem}
So, if two varieties belong to the same Hilbert scheme ${\rm
Hilb}^{P(x)}_{\mathbf P^N}$, then they appear as fibers of the
universal family $u : \mathcal U \longrightarrow {\rm
Hilb}^{P(x)}_{\mathbf P^N}$, in which ${\rm Hilb}^{P(x)}_{\mathbf
P^N}$ is {\it connected}. In this way, they are connected by
projective flat deformation.

Let $Z$ be a Calabi--Yau threefold and let $H$ be an ample divisor
on $Z$. Then, by the Kodaira vanishing theorem and the
Riemann-Roch formula, we have
$$\dim H^0(\mathcal O_Z(nH)) = \chi(\mathcal O_Z(nH)) =
\frac{H^3}{6} n^3 + \frac{H \cdot c_2(Z)}{12}n\,\, .$$ Here
$c_2(Z) = c_2(T_Z)$ is the second Chern class of $Z$. It is also
known that $10H$ is always very ample on $Z$ (\cite{OP}).
Therefore, as a special case of Theorem 0.2, one obtains the
following:
\begin{theorem}\label{key} Two Calabi--Yau threefolds have
the same Hilbert polynomial, belong to the same Hilbert scheme of
some projective space, and accordingly connected by projective
flat deformation, if and only if they have ample divisors that have
the same values of \cen{$H^3$ and $H \cdot c_2$.}
\end{theorem}

In general, two Calabi--Yau threefolds are unlikely to be connected
by projective flat deformation, especially if they are of
different topological type. Let $X$ and $Y$ be a complete
intersection of two cubics in $\mathbf P^5$ and a quintic
hypersurface  in $\mathbf P^4$ respectively. Then we always have
$$9 k^3 =(kH_X)^3 \neq (lH_Y)^3 = 5 l^3$$
for any positive integers $k, l$, where $H_X$ and $H_Y$ are the
ample generators of the Picard groups of $X$ and $Y$ respectively.
So $X$ and $Y$ can not be connected by any projective flat
deformation.

Our Calabi--Yau threefolds in Theorem \ref{mainthm} are the
famous rigid Calabi--Yau threefold $X_{\phi}$
constructed by Beauville \cite{Be2} and its birational
modification $X_{T}$ studied by the second named author \cite{Og}
(See also Section \ref{beauville}).

The structure of this note is as follows: We discuss some toy case
of elliptic curves in Section \ref{toy}. This explains some idea
behind our consideration. In Section \ref{beauville}, we recall
Beauville's rigid Calabi--Yau threefold $X_{\phi}$ and its
birational modification $X_{T}$. Sections \ref{topology} and
\ref{main} are devoted to the proof of Theorem \ref{mainthm}.

\par \vskip 1pc

{\it Acknowledgement.} We would like to express our thanks to
Professors J.M. Hwang, J.H. Keum, B. Kim for valuable discussions.
\par \vskip 1pc

\section{Toy example: connecting two elliptic curves in two ways} \label{toy}
Let $C_{\lambda}$ ($\lambda \neq 0, 1$) be the elliptic curve defined by the Weierstrass equation
$$y^2  = x(x-1)(x-\lambda)\,\, .$$
Obviously, any two elliptic curves $C_{\lambda_1}$ and $C_{\lambda_2}$
are connected by the following projective smooth family:
$$\psi : \mathcal X = \{([x_0 : x_1 : x_2], \lambda) \in \mathbf P^2
 \times \mathcal B\, {\big |}\, x_1^2x_2  - x_0(x_0-x_2)(x_0- \lambda x_2) =0\}
 \longrightarrow \mathcal B\,\, .$$
Here and hereafter, we put $\mathcal B = \mathbf P^1 \setminus
\{0, 1, \infty\}$.

Yet, we can connect  $C_{\lambda_1}$ and  $C_{\lambda_2}$ by another way.

 Let $D$ be a hyperelliptic curve with a hyperelliptic involution $\iota$
 and let $\Xi$ be the set of the branch points of $\iota$
 in $D/\langle \iota \rangle \simeq \mathbf P^1$. We consider
 the natural morphisms,
$$\varphi_1: \widetilde{C_{\lambda_1} \times  D \big{/}
\langle(-1, \iota) \rangle }  \longrightarrow D/\langle \iota
\rangle \simeq \mathbf P^1$$ and
$$\varphi_2: \widetilde{C_{\lambda_2} \times  D \big{/}
\langle (-1, \iota) \rangle} \longrightarrow D/\langle \iota
\rangle \simeq \mathbf P^1.$$ Here $\widetilde{\hspace{.5cm}}$'s
are the minimal resolutions. We regard $\varphi_1$ and $\varphi_2$
as projective flat deformations. Then, for $q \in \Xi$, the
scheme-theoretic fiber $\varphi_1^{-1}(q)= 2l + l_0 + l_1 +
l_\infty + l_{\lambda_1}$ consists of 5 $\mathbf P^1$'s,
intersecting like:
$$ 
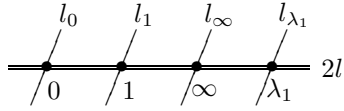
\xy
(5,-.2)*{}; (45,-.2)*{} **\dir{-};
(5,.2)*{}; (45,.2)*{} **\dir{-};
(8,-5)*{}; (12,5)*{} **\dir{-};
(18,-5)*{}; (22, 5)*{} **\dir{-};
(28,-5)*{}; (32,5)*{} **\dir{-};
(38,-5)*{}; (42, 5)*{} **\dir{-};
(10,0)*{\bullet}="A"+(1,-3)*{0}+(2, 10)*{l_0};
(20,0)*{\bullet}="B"+(1,-3)*{1}+(2, 10)*{l_1};
(30,0)*{\bullet}="C"+(1,-3)*{\infty}+(2, 10)*{l_\infty};
(40,0)*{\bullet}="D"+(1,-3)*{\lambda_1}+(2, 10)*{l_{\lambda_1}}+(5,-7)*{2l};
\endxy $$
and $\varphi_1^{-1}(p) \simeq C_{\lambda_1}$ for $p \notin \Xi$.
Similarly the scheme-theoretic fiber $\varphi_2^{-1}(q)$ for $q \in \Xi$ is like:
$$ 
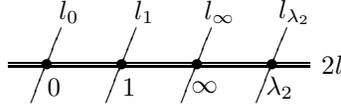
\xy
(5,-.2)*{}; (45,-.2)*{} **\dir{-};
(5,.2)*{}; (45,.2)*{} **\dir{-};
(8,-5)*{}; (12,5)*{} **\dir{-};
(18,-5)*{}; (22, 5)*{} **\dir{-};
(28,-5)*{}; (32,5)*{} **\dir{-};
(38,-5)*{}; (42, 5)*{} **\dir{-};
(10,0)*{\bullet}="A"+(1,-3)*{0}+(2, 10)*{l_0};
(20,0)*{\bullet}="B"+(1,-3)*{1}+(2, 10)*{l_1};
(30,0)*{\bullet}="C"+(1,-3)*{\infty}+(2, 10)*{l_\infty};
(40,0)*{\bullet}="D"+(1,-3)*{\lambda_2}+(2, 10)*{l_{\lambda_2}}+(5,-7)*{2l};
\endxy $$
and $\varphi_2^{-1}(p) \simeq C_{\lambda_2}$ for $p \notin \Xi$.
The singular schemes $\varphi_1^{-1}(q)$ and $\varphi_2^{-1}(q)$
can be put into a projective flat family, in which the fibers are
of the form:
$$ 
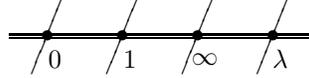
\xy
(5,-.2)*{}; (45,-.2)*{} **\dir{-};
(5,.2)*{}; (45,.2)*{} **\dir{-};
(8,-5)*{}; (12,5)*{} **\dir{-};
(18,-5)*{}; (22, 5)*{} **\dir{-};
(28,-5)*{}; (32,5)*{} **\dir{-};
(38,-5)*{}; (42, 5)*{} **\dir{-};
(10,0)*{\bullet}="A"+(1,-3)*{0};
(20,0)*{\bullet}="B"+(1,-3)*{1};
(30,0)*{\bullet}="C"+(1,-3)*{\infty};
(40,0)*{\bullet}="D"+(1,-3)*{\lambda};
\endxy $$
For example, the natural projection $\psi: \mathcal Y \longrightarrow \mathcal B$, where
$$\mathcal Y = \{ ([x_0 : x_1], [y_0: y_1], \lambda) \in
\mathbf P^1 \times \mathbf P^1 \times \mathcal B\,  {\big |}\,
x_0^2 y_0 y_1 (y_0 - y_1) (y_0-\lambda_1 y_1) =0 \}$$
is such a family. In this way, $C_{\lambda_1}$ and
 $C_{\lambda_2}$ are connected by a chain of three projective flat deformations.

In the second method, smooth fibers in families are only
$C_{\lambda_1}$ and $C_{\lambda_2}$ and they are connected through
{\it very singular spaces}. So, the method suggests some
possibilities to connect two rigid manifolds of different
topological structure. This is the idea behind our construction.

\section{Beauville's rigid Calabi--Yau threefold and its modification} \label{beauville}
We briefly recall the two rigid Calabi--Yau threefolds $X$ and $Y$ that appear
in Theorem \ref{mainthm}.

Let $\zeta = e^{2\pi \sqrt{-1}/3}$. By $E_\zeta$, we denote the elliptic curve whose
period is $\zeta$ and by $E_\zeta^n/\langle \zeta \rangle $ the quotient
$n$-fold of
the product manifold $E_\zeta^n$ by the scalar multiplication by
$\zeta$. Let \cen{$Q_0 =0$, $Q_1 = (1-\zeta)/3$ and $Q_2
=-(1-\zeta)/{3}$} in $E_\zeta$. These are exactly the fixed points
of the scalar multiplication by $\zeta$ on $E_\zeta$. For $i_k=0, 1, 2$, let
$$Q_{i_1 i_2 \cdots i_n} = (Q_{i_1}, Q_{i_2}, \cdots, Q_{i_n}) \in E_\zeta^n$$
and let $\overline Q_{i_1 i_2 \cdots i_n}$ be its image in
$E_\zeta^n/\langle \zeta \rangle $. Then $\overline X =
E_\zeta^3/\langle \zeta \rangle $ has singularities of type
$\frac{1}{3}(1,1,1)$ at $\overline Q_{ijk}$'s and  the blow-up
$\pi : X_{\varphi} \longrightarrow \overline X$ at these 27
singular points gives a Calabi--Yau threefold $X_\phi$. This is the
famous rigid Calabi--Yau threefold found by Beauville \cite{Be2}.
We denote by $E_{ijk}$ the exceptional divisor lying over
$\overline Q_{ijk}$. The surfaces $E_{ijk}$ is isomorphic to
$\mathbf P^{2}$.

\begin{figure}[h]
$$ \xy
(0,0)*{}="A";
(50,0)*{}="B";
(50,50)*{}="G";
(0,50)*{}="Y";
"A";"B" **\dir{-};
"Y";"G" **\dir{-};
"Y";"A" **\dir{-};
"B";"G" **\dir{-};
(15,3)*{}; (15,47)*{} **\dir{-};
(15,10)*{\bullet}="A"+(9,-3)*{E_{ij0}}+(20, 20)*{E_\zeta} ;
(15,25)*{\bullet}="B"+(9,-3)*{E_{ij1}}+(34, 3)*{X_\phi};
(15,40)*{\bullet}="C"+(9,-3)*{E_{ij2}}+(-6,10)*{l_{ij}};
(5,5)*{}="a";
(17,5)*{}="b";
(13,15)*{}="d";
(25,15)*{}="c";
"a";"b" **\dir{-};
"b";"c" **\dir{-};
"c";"d" **\dir{-};
"d";"a" **\dir{-};
(5,20)*{}="a";
(17,20)*{}="b";
(13,30)*{}="d";
(25,30)*{}="c";
"a";"b" **\dir{-};
"b";"c" **\dir{-};
"c";"d" **\dir{-};
"d";"a" **\dir{-};
(5,35)*{}="a";
(17,35)*{}="b";
(13,45)*{}="d";
(25,45)*{}="c";
"a";"b" **\dir{-};
"b";"c" **\dir{-};
"c";"d" **\dir{-};
"d";"a" **\dir{-};
(40,0)*{}; (40,50)*{} **\dir{-};
(15,-15)*{\ast}="A"+(0,-3)*{\overline Q_{ij}}+(45,3)*{B};
(40,-15)*{\ast}="B"+(0,-3)*{P}+(0, -5)*{(P \notin \bigcup_{i,j}\{ \overline Q_{ij}\} )};
(0,-15)*{}; (50,-15)*{} **\dir{-};
{\ar(25,-3)*{}; (25,-13)*{} }+(28,6)*{p_\phi};
\endxy $$
\caption{}\label{X_phi}
\end{figure}

Let
$$p_\phi : X_\phi \longrightarrow B := E_\zeta^2/\langle \zeta \rangle$$
be the morphism, induced by the projection ${\rm pr}_{12} :
E_\zeta^3\longrightarrow E_\zeta^2$. Then we have
$$p_{\phi}^{-1}{\overline Q_{ij}} = l_{ij} \cup E_{ij0}\cup E_{ij1}\cup E_{ij3}\,\, .$$
Here $l_{ij}$ is a smooth rational curve meeting $E_{ijk}$
transversally. See Figure \ref{X_phi}. The normal bundle of
$l_{ij}$ in $X_\phi$ is:
$$N_{X_\phi|l_{ij}} = \mathcal O_{l_{ij}}(-1)^{\oplus 2}.$$
Performing the elementary transformation along $\bigcup_{i,j}
l_{ij}$, we obtain a smooth threefold $X_T$. This $X_T$
corresponds to that in \cite{Og} for $T=\{(i,j) | i, j=0, 1, 2\}$.
Denote the proper transform of $E_{ijk}$ in $X_T$ by $F_{ijk}$.
Note that $F_{ijk}$ is the first Hirzebruch surface $\mathbf F_1$.
Compare Figure \ref{X_T} with Figure \ref{X_phi}.

Now we summarize some properties of $X_\phi$ and $X_T$, showed in
\cite{Be2} and \cite{Og}:
\begin{theorem}

\begin{enumerate}
\item $X_\phi$ and $X_T$ are both projective and simply-connected.
Accordingly they are in fact Calabi--Yau threefolds. \item
$h^{1,2}(X_\phi)= h^{1,2}(X_T) = 0$, i.e. $X_\phi$ and $X_T$ are
rigid.
\end{enumerate}
\end{theorem}

So, $X_\phi$ and $X_T$ are birational, rigid Calabi--Yau
threefolds. In fact, these $X_\phi$ and $X_T$ are the Calabi--Yau
threefolds $X$ and $Y$ in our Theorem \ref{mainthm}. We shall show
that $X_\phi$ and $X_T$ are non-homeomorphic in Section
\ref{topology} and that $X_\phi$ and $X_T$ are connected by
projective flat deformation in Section \ref{main}.

\begin{figure}[h]
$$ \xy
(-10,-10)*{}="A";
(-10,-6)*{}="A'";
(-10, 6)*{}="A''";
(-10, 10)*{}="B";
(10, 10)*{}="C";
(10, -10)*{}="D"+(1, -3)*{F_{ij0}};
"A";"A'" **\dir{-};
"A'";"A''" **\dir{..};
"A''";"B" **\dir{-};
"B";"C" **\dir{-};
"C";"D" **\dir{-};
"D";"A" **\dir{-};
(-16, -6)*{}="A";
(-10, 0)*{}="A'";
(-5, 5)*{}="B";
(10, 5)*{}="B'";
(15, 5)*{}="C"+(2,2)*{F_{ij1}};
(4, -6)*{}="D";
"A";"A'" **\dir{-};
"A'";"B" **\dir{..};
"B";"B'" **\dir{..};
"B'";"C" **\dir{-};
"C";"D" **\dir{-};
"D";"A" **\dir{-};
(-5, -5)*{}="A";
(-10, 0)*{}="A'";
(-16, 6)*{}="B";
(4, 6)*{}="C";
(15, -5)*{}="D"+(4,-2)*{F_{ij2}}+(12,6)*{X_T};
(10, -5)*{}="D'";
"A";"A'" **\dir{..};
"A'";"B" **\dir{-};
"B";"C" **\dir{-};
"C";"D" **\dir{-};
"D";"D'" **\dir{-};
"D'";"A" **\dir{..};
(-10, 0)*{}; (10, 0)*{} **\dir{-};
(-25,-25)*{}="A";
(25,-25)*{}="B";
(25,25)*{}="G";
(-25,25)*{}="Y";
"A";"B" **\dir{-};
"Y";"G" **\dir{-};
"Y";"A" **\dir{-};
"B";"G" **\dir{-};
(23,-25)*{}; (23,25)*{} **\dir{-};
(0,-40)*{\ast}="A"+(0,-3)*{\overline Q_{ij}}+(35,3)*{B };
(23,-40)*{\ast}="B"+(0,-3)*{P}+(0, -5)*{(P \notin \bigcup_{i,j} \{\overline Q_{ij}\} )};
(-25,-40)*{}; (25,-40)*{} **\dir{-};
{\ar(0,-28)*{}; (0,-38)*{} }+(28,6)*{p_T};
(30,20)="R"+(1,3)*{E_\zeta};
{\ar@/^.3pc/"R";(23.7,10)};
\endxy $$
\caption{}\label{X_T}
\end{figure}

Here, we summarize notations which will be frequently used in the next two sections:

\begin{notation}\label{lst}

\begin{itemize}
\item $\zeta = e^{2\pi \sqrt{{-}1}/3}$, the primitive third root
of unity in the upper half plane. \item $E_\zeta = \mathbf C /
(\mathbf Z \oplus \mathbf Z \zeta)$ is the elliptic curve with
period $\zeta$. \item $Q_0 =0$, $Q_1 = (1-\zeta)/3$ and $Q_2
=-(1-\zeta)/{3}$ in $E_\zeta$: the fixed points of the scalar
multiplication by $\zeta$ on $E_\zeta$. \item $Q_{i_1 i_2 \cdots
i_n} = (Q_{i_1}, Q_{i_2}, \cdots, Q_{i_n}) \in E_\zeta^n$ for $
i_1, i_2, \cdots, i_n \in \{ 0, 1, 2 \}$. \item $\overline Q_{i_1
i_2 \cdots i_n}$ is the image of $Q_{i_1 i_2 \cdots i_n}$ in
$E_\zeta^n/\langle \zeta \rangle$. \item $\overline X =
E_\zeta^3/\langle \zeta \rangle $, $B = E_\zeta^2/\langle \zeta
\rangle $. Here $E_\zeta^2$ in the definition of $B$ is the
product of the first two factors of $E_\zeta^3$. \item
$q:E_\zeta^3 \longrightarrow \overline X$ is the quotient map.
\item ${\rm pr}_{i}:E_\zeta^3 \longrightarrow E_\zeta$ is the
projection to the $i^{\rm th}$ factor. \item ${\rm
pr}_{ij}:E_\zeta^3 \longrightarrow E_\zeta^2$ is the projection to
the product of $i^{\rm th}$ and $j^{\rm th}$ factors. \item $
p_{ij}:\overline X \longrightarrow E_\zeta^2 / \langle \zeta
\rangle$ and  $p_{i}:\overline X \longrightarrow E_\zeta / \langle
\zeta \rangle$ are the the morphisms induced by ${\rm pr}_{ij}$
and ${\rm pr}_{i}$ respectively. \item $g_i : B \longrightarrow
E_\zeta/\langle \zeta \rangle = \mathbf P^1$ is the morphism,
induced by the projection $E_\zeta^2 \longrightarrow E_\zeta$ to
the $i^{\rm th}$ factor ($i = 1$, $2$). \item $\pi:X_\phi
\longrightarrow \overline X$ is the blow-up at $\{ \overline
Q_{ijk}| i,j,k = 0, 1, 2\}$. \item $p_\phi =   p_{12} \circ \pi
:X_\phi \longrightarrow B$. \item $p_T :X_T \longrightarrow B$ is
the projection, induced by $p_\phi$. \item $E_{ijk} \simeq \mathbf
P^2$ is the exceptional divisor over $\overline Q_{ijk}$ by the
blow-up $\pi : X_\phi \longrightarrow \overline X$. \item $F_{ijk}
\simeq \mathbf F_1$ is the proper transformation of $E_{ijk}$ in
$X_T$.
\end{itemize}
\end{notation}

The next lemma will be also frequently used in the next two sections:
\begin{lemma}\label{chn} Let $Z$ be a Calabi--Yau threefold and let $D$ be a smooth
divisor on $Z$. Then, $D^3 = c_1(T_D)^2$ and $D \cdot c_2(Z) =
-c_1(T_D)^2 + c_2(T_D)$.
\end{lemma}
We also note that $c_1(T_D)^2 = K_D^{2}$ and that $c_2(T_D) =
c_2(D)$ is the topological Euler number of the surface $D$.

\begin{proof} This follows from the fact that $c_{1}(Z) = 0$ and the normal sequence
$$0 \longrightarrow T_D \longrightarrow T_{Z}\vert D \longrightarrow N_{Z|D}
\longrightarrow 0\,\, .$$
\end{proof}

\section{Topological difference between $X_\phi$ and $X_T$}\label{topology}
In this section we shall prove (\ref{non-homeo}) of Theorem
\ref{mainthm}, i.e. that $X_\phi$ and $X_T$ are not homeomorphic.
Since the linear form $c_{2}(Z) : H^{2}(Z, \mathbf Z)
\longrightarrow \mathbf Z$ and the cubic form $c_{Z} : {\rm
Sym}^{3} H^{2}(Z, \mathbf Z) \longrightarrow \mathbf Z$ are
topological invariants, the result follows from:

\begin{theorem}\label{cbc}

(1) The linear form given by $c_2(X_{\phi})$ is divisible by $6$, i.e.
$D \cdot c_2(X_{\phi}) \equiv 0$ (${\rm mod}\, 6$) for each $D \in H^2(X_{\phi}, \mathbf Z)$, while the linear form $c_2(X_{T})$ is not.

(2) The cubic form of $X_{\phi}$ is divisible by $3$, i.e. $D^{3} \equiv 0$ (${\rm mod}\, 3$) for each $D \in H^2(X_{\phi}, \mathbf Z)$, while the cubic form of $X_{T}$ is not.

\end{theorem}

\begin{remark} As far as we know, Friedman is the first who found
a pair of birational projective Calabi--Yau threefolds which are
not homeomorphic (\cite{Fr} Example 7.7.). His examples are based
on [Sc] and they are not rigid. Our proof here is inspired by his
argument there.
\end{remark}

We shall prove Theorem \ref{cbc} in the sequel.

Let $F \simeq \mathbf F_1$ be one of $F_{ijk}$ in $X_T$. Then, by Lemma
\ref{chn},
$$F^3 = K_{F}^{2} = 8\,\, ,\,\, F \cdot c_2(X_T) = c_2(F) - K_{F}^{2} = 4 - 8 = -4\, .$$
Clearly, none of them is divisible by $3$.

In the rest of this section, we shall show $6$-divisibility of the linear
form $c_{2}(X_\phi)$. $3$-divisibility of the cubic form then follows from the
Riemann-Roch formula (cf. Introduction). Here we note that
${\rm Pic}\, X_{\phi} \simeq H^{2}(X_{\phi}, \mathbf Z)$. {\it From now
until the end of this section, we write}
$$E = E_{\zeta}\,\, , \,\, X = X_{\phi}\,\,.$$
For other notations, see Notation \ref{lst}.
\begin{proposition}\label{gen} Let $n \ge 2$.

(1) The N\'eron-Severi group ${\rm NS}\, (E^2)$ is generated by
the classes of the four divisors, $\{0\} \times E$, $E \times\{0\}$, $\Delta$,
$\Gamma$. Here $\Delta$ is the diagonal and $\Gamma$ is the graph of
the automorphism $\zeta : E \longrightarrow E$.

(2) The N\'eron-Severi group ${\rm NS}\, (E^n)$ is generated by the subgroups
${\rm pr}_{ij}^{*}NS(E^2)$ ($1 \le i < j \le n$).
\end{proposition}

\begin{proof} The four
classes in (1) are clearly in ${\rm NS}(E^2)$, and their intersection matrix is

$$\left(\begin{array}{rrrr} 0 & 1& 1 & 1\\
1 & 0& 1 & 1\\
1 & 1& 0 & 3\\
1 & 1& 3 & 0
\end{array} \right)\,\, .$$

The discriminant of this matrix is $3$. On the other hand, the discriminant
of the transcendental lattice of $E^{2}$ is $3$ by \cite{ShM}. Thus,
the discriminant of ${\rm NS}(E^2)$ is also $3$. Since ${\rm NS}(E^2)$
is torsion free, the assertion (1) follows.

Let us show (2). By the K\"unneth formula, we have
$$H^2(E^n, \mathbf Z) = \bigoplus_{i=1}^{n} {\rm pr}_{i}^{*}H^2(E,
\mathbf Z) \oplus \bigoplus_{1 \le i < j \le n} {\rm pr}_{ij}^{*}
(H^1(E, \mathbf Z) \otimes H^{1}(E, \mathbf Z))\,\, .$$
This decomposition is compatible with the Hodge decomposition.
Since ${\rm NS}(E^n) = H^2(E^n, \mathbf Z) \cap H^{1,1}(E^n)$ by the
Lefschetz $(1,1)$-theorem,
we have then
$${\rm NS}\, (E^n) = \bigoplus_{i=1}^{n} {\rm pr}_{i}^{*}H^2(E, \mathbf Z)
\oplus \bigoplus_{1 \le i < j \le n} {\rm pr}_{ij}^{*}(H^1(E, \mathbf Z)
\otimes H^{1}(E, \mathbf Z) \cap H^{1,1}(E^2))\,\, .$$
Again, by the Lefschetz $(1,1)$-theorem, the groups ${\rm
pr}_{k}^{*}H^{2}(E, \mathbf Z)$ ($k =i$, $j$) and ${\rm
pr}_{ij}^{*}(H^1(E, \mathbf Z) \otimes H^{1}(E, \mathbf Z) \cap
H^{1,1}(E^2))$ are subgroups of ${\rm pr}_{ij}^{*}{\rm NS}(E^2)$,
in which $E^{2}$ is the product of $i^{\rm th}$ and $j^{\rm th}$
factors of $E^n$. This implies (2).
\end{proof}

Recall that $\overline{X} := E^3/\langle \zeta \rangle$. In
particular, $\overline{X}$ is $\mathbf Q$-factorial. A bit more
precisely, the divisor $3D$ is Cartier for any Weil divisor $D$ on
$\overline{X}$. Let $N^{1}(\overline{X})$ be the group of the
numerically equivalent classes of Weil divisors on $\overline{X}$.
Note that Cartier divisors and Weil divisors are the same on $E^3$
or on $X$ (as $E^3$ and $X$ are smooth) and the numerical
equivalence and the algebraic equivalence of divisors are also the
same on $E^3$ or on $X$ (as their N\'eron-Severi groups are
torsion free).

\begin{proposition}\label{ism} The group homomorphism $q^* : N^1(\overline{X})
\longrightarrow {\rm NS}\, (E^3)$ is an isomorphism.
\end{proposition}

\begin{proof} Our argument here is similar to \cite{Nm}. Since $\overline{X}$
is $\mathbf Q$-factorial and $q$ is finite, the group homomorphism $q^{*}$ is indeed
well-defined and injective.

Let $[H] \in {\rm NS}(E^3)$. We need to find $D \in N^1(\overline{X})$ such that $[H] = q^{*}D$.

{\it Claim 1.} We can (and will) choose the representative $H \in
{\rm Pic}\, E^3$ of the class $[H]$, such that $\zeta^{*}H = H$ as
line bundles.

{\it Proof of Claim 1.} Take the origin of $E$ as polarization of
$E$. One can then identify ${\rm Pic}^{0}(E) = E$ in an
equivariant way with respect to the action of $\zeta$. Under this
identification, we have an identification ${\rm Pic}^{0}(E^3) =
E^3$ in which the action of $\zeta^{*}$ on ${\rm Pic}^{0}(E^3)$ is
the same as the diagonal action $(\zeta, \zeta, \zeta)$ on $E^3$.
Note also that $\zeta^{*} = {\rm {id}}$ on ${\rm NS}(E^3)$ as $\zeta^{*} =
{\rm{id}}$ on the wider space $H^{1,1}(E^3)$.

Put $T = \zeta^{*} H - H$. Here the equality is as line bundles.
Then $T = (T_1, T_2, T_3)$ is an element of ${\rm Pic}^0(E^3) =
E^3$, as $\zeta^{*}[H] = [H]$. Note that there is a point $P =
(P_1, P_2, P_3) \in {\rm Pic}^{0}(E^3) = E^3$ such that
$$(P_{1}, P_{2}, P_{3}) - (\zeta P_1, \zeta P_2, \zeta P_3) = (T_1, T_2, T_3)\,\, .$$
The line bundle $H + P$ is a desired representative. q.e.d. for
Claim 1.

From now, we regard $H$ as an effective divisor in $\vert H \vert$ rather than the line bundle.

{\it Claim 2.} We may (and will) assume that there is an effective
divisor $H$ in $\vert H \vert$ such that $\zeta^{*}H = H$ as
divisors.

{\it Proof of Claim 2.} Since $q$ is finite, the divisor $q^{*}A$
is ample if $A$ is ample. Thus, by adding $q^{*}A$ with
sufficiently ample $A$ to $H$, we may assume that $\vert H \vert$
is a free linear system. Since $\zeta^{*}H = H$ as line bundles,
$\zeta$ acts on the projective space $\vert H \vert$. This action
certainly has a fixed points. Let $H$ be a divisor corresponding
to (one of) the fixed point. Then $\zeta^{*} H = H$ as divisors on
$E^3$. q.e.d. for Claim 2.

Let $\overline{H} = q_{*}H$ as Weil divisors. Since $\zeta^{*} H = H$ as divisors
and $(E^3)^{\langle \zeta \rangle}$ consists of finitely many points, there
is a divisor $D$ such that
$\overline{H} = 3D$ as Weil divisors. For this $D$, we have
$$3 q^* D = q^{*}\overline{H} = H + \zeta^{*} H + (\zeta^{*})^{2} H = 3H\,\, .$$
Since ${\rm NS}\, (E^3)$ is torsion free, this implies $q^* D = H$.
\end{proof}

\begin{proposition}\label{wel} Let $\tilde{D}_{ijl}$ ($1 \le l \le 4$)
be the divisors on $E^{3}$, which are pull back of the four divisors
$E \times \{0\}$, $\{0 \} \times E$, $\Delta$, $\Gamma$ on $E^2$ by $p_{ij}$
($1 \le i < j \le 3$).
Let $\overline{D}_{ijl} := (q_{*}\tilde{D}_{ijl})_{\rm red}$. Then, the (classes of)
$12$ Weil divisors $\overline{D}_{ijl}$ generate $N^1(\overline{X})$.
\end{proposition}

\begin{proof} We note that $\zeta^{*}\tilde{D}_{ijl} = \tilde{D}_{ijl}$
as divisors on $E^3$. Thus $\tilde{D}_{ijl} = \pi^{*}\overline{D}_{ijl}$
(cf. Proof of Proposition \ref{ism}). Since $\tilde{D}_{ijl}$ generate
${\rm NS}(E^3)$ by Proposition \ref{gen}, the result follows from Proposition \ref{ism}.
\end{proof}

Let $D_{ijl}$ be the proper transform of $\overline{D}_{ijl}$ on
$X$ by $\pi : X \longrightarrow \overline X$.

\begin{proposition}\label{fidx} ${\rm NS}\, (X)$ is contained in
the subgroup of ${\rm NS}\, (X) \otimes \mathbf Q$ generated by the classes
of the following divisors:
$$D_{ijl}\,\, ,\,\, E_{ijk}\,\, , \,\,
T_{\Lambda, \Lambda'} := \frac{1}{3}\sum_{(i, j, k) \in \Lambda}
E_{ijk} + \frac{1}{3}\sum_{(i', j', k') \in \Lambda'}
2E_{i'j'k'}\,\, ,$$ where $\Lambda$ and $\Lambda'$ are some disjoint
subsets (possibly empty) of the product set $\{0, 1, 2\}^{3}$ such
that $\Lambda \cap \Lambda' = \emptyset$ and such that
both $\vert \Lambda
\vert$ and $\vert \Lambda' \vert$ are divisible by $3$.
\end{proposition}

\begin{proof}

Let $D$ be a prime divisor on $X$. Put $\overline{D} = \pi_{*}D$
as Weil divisors. Then, by Proposition \ref{wel}, there are
integers $b_{ijl}$ such that $\overline{D} = \sum_{i,j, l}
b_{ijl}\overline{D}_{ijl}$ in $N^{1}(\overline{X})$. Since
$3\overline{D}$ and $3\overline{D}_{ijl}$ are Cartier, there are
integers $a_{ijk}$ such that
$$D = \sum_{i,j, l} b_{ijl}D_{ijl} + \frac{1}{3} \sum_{i, j, k}
a_{ijk}E_{ijk}\,\,$$
in ${\rm NS}(X) \simeq {\rm Pic}\, X$. So, the result follows from the
next lemma. \end{proof}

\begin{lemma}\label{mult} Let $\Lambda$ and $\Lambda'$ be subset of $\{0, 1, 2\}^{3}$
such that $\Lambda \cap \Lambda' = \emptyset$. If
$$M := \sum_{(i, j, k) \in \Lambda} E_{ijk} + \sum_{(i', j', k') \in \Lambda'}
2E_{i'j'k'}$$ is divisible by $3$ in ${\rm Pic}\, X$, then both
$\vert \Lambda \vert$ and $\vert \Lambda' \vert$ are divisible by
$3$.
\end{lemma}

\begin{proof}
Let $\alpha \in \{0, 1, 2\}$. Let $D_{\alpha}$ be the divisor on
$X$, which is the proper transform of the divisor
$\overline{D}_{\alpha} = (p_{3}^{*}Q_{\alpha})_{\rm red}$ on
$\overline{X}$ (See Figure \ref{o2} and Notation \ref{lst} for
$Q_{\alpha}$).
\begin{figure}[h]
$$
\xy
(0,0)="A";
(20,-20)="B";
(20,0)="C";
(0,20)="D"+(-2,3)*{\overline D_1};
"A";"B" **\dir{-};
"B";"C" **\dir{-};
"C";"D" **\dir{-};
"D";"A" **\dir{-};
(2,3)*{\bullet};
(2,8)*{\bullet};
(2,13)*{\bullet};
(18,-15)*{\bullet};
(18,-10)*{\bullet};
(18,-5 )*{\bullet};
(10,-6)*{\bullet};
(10,-1)*{\bullet};
(10,4 )*{\bullet};
(20,-3)="A";
(15,5)="A'";
(35,-18)="B";
(35,2)="C";
(15,22)="D"+(-2,3)*{\overline D_2};
"A";"B" **\dir{-};
"B";"C" **\dir{-};
"C";"D" **\dir{-};
"D";"A'" **\dir{-};
(17,5)*{\bullet};
(17,10)*{\bullet};
(17,15)*{\bullet};
(33,-13)*{\bullet};
(33,-8)*{\bullet};
(33,-3 )*{\bullet};
(25,-4)*{\bullet};
(25,1)*{\bullet};
(25,6 )*{\bullet};
(35,-1)="A";
(30,7)="A'";
(50,-16)="B";
(50,4)="C";
(30,24)="D"+(-2,3)*{\overline D_3};
"A";"B" **\dir{-};
"B";"C" **\dir{-};
"C";"D" **\dir{-};
"D";"A'" **\dir{-};
(32,7)*{\bullet};
(32,12)*{\bullet};
(32,17)*{\bullet};
(48,-11)*{\bullet};
(48,-6)*{\bullet};
(48,-1 )*{\bullet};
(40,-2)*{\bullet};
(40,3)*{\bullet};
(40,8 )*{\bullet};
\endxy
$$\caption{}\label{o2}
\end{figure}

Since $\overline{D}_{\alpha}$ passes through $9$-singular points of $\overline{X}$,
the surface $D_{\alpha}$ meets the $9$-exceptional divisors, say,
$$E_{00\alpha}\, ,\, E_{01\alpha}\, ,\, E_{02\alpha}\, ,\, E_{10\alpha}\, ,
\, E_{11\alpha}\, ,\,E_{12\alpha}\, ,\, E_{20\alpha}\, , \,
E_{21\alpha}\, , \, E_{22\alpha}\, .$$ We put $l_{ij \alpha} :=
E_{ij \alpha} \vert_{D_{\alpha}}$. These are all $(-3)$-curves.
The surface $D_{\alpha}$ is a non-relatively minimal rational
elliptic surface with $3$-singular fibers $l_{i0\alpha} +
l_{i1\alpha} + l_{i2\alpha} + 3c_{i}$ ($i = 0$, $1$, $2$) as in
the figure below (Figure \ref{o1}). Here $c_{i}$ are
$(-1)$-curves. Let $\nu_{\alpha} : D_{\alpha} \longrightarrow
D_{\alpha}'$ be the contraction of the three $(-1)$-curves
$c_{i}$. Let $l_{ij \alpha}' = \nu(l_{ij \alpha})$. Then
$D_{\alpha}'$ is a relatively minimal rational elliptic surface
with $3$ singular fibers $l_{i0\alpha}' + l_{i1\alpha}' +
l_{i2\alpha}'$ (Figure \ref{o1}).
\begin{figure}[h]
$$
\xy
(-10,-10)*{}="A";
(10, -10)*{}="B";
(10, 10)*{}="C";
(-10,10)*{}="D";
"A";"B" **\dir{-};
"B";"C" **\dir{-};
"C";"D"**\dir{-};
"D";"A" **\dir{-};
(0,8)="E"+(4.5,6)*{c_1};
(0,-8)="F";
"E";"F" **\dir{-};
(-2,4);(2,6) **\dir{-};
(-2,-1);(2,1)**\dir{-};
(-2,-6);(2,-4) **\dir{-};
(-6,8)="E"+(4.5,6)*{c_0};
(-6,-8) ="F"+(-8,13)*{D_\alpha};
"E";"F" **\dir{-};
(-8,4);(-4,6)**\dir{-};
(-8,-1);(-4,1) **\dir{-};
(-8,-6);(-4,-4) **\dir{-};
(6,8)="E"+(4.5,6)*{c_2};
(6,-8) ="F";
"E";"F" **\dir{-};
(4,4);(8,6) **\dir{-};
(4,-1);(8,1) **\dir{-}; (4,-6);(8,-4)
**\dir{-};
{\ar@/_.1pc/(3, 13);(0.5,8)}; {\ar@/_.1pc/(9, 13);(6.5,8)};
{\ar@/_.1pc/(-3, 13);(-5.5,8)}; {\ar (20,0)*{}; (30,0)*{}};
(5,-15)="A"+(1,-3)*{l_{1 j \alpha}}; {\ar@/_.1pc/"A"; (1,-5)*{}};
(40,-10)*{}="A";
(60, -10)*{}="B";
(60, 10)*{}="C";
(40, 10)*{}="D";
"A";"B" **\dir{-};
"B";"C" **\dir{-};
"C";"D" **\dir{-};
"D";"A" **\dir{-};
(50,6)="E";
(50,-6) ="F";
"E";"F" **\dir{-};
(48,5);(52,-5) **\dir{-};
(52,5);(48,-5) **\dir{-};
(44,6)="E";
(44,-6) ="F"+(-8,13)*{D_\alpha'};
"E";"F" **\dir{-};
(42,5);(46,-5) **\dir{-};
(46,5);(42,-5) **\dir{-};
(56,6)="E";
(56,-6) ="F";
"E";"F" **\dir{-};
(54,5);(58,-5) **\dir{-};
(58,5);(54,-5) **\dir{-};
(55,-15)="A"+(1,-3)*{l_{1 j \alpha}'}; {\ar@/_.1pc/"A";
(50.5,-5)*{}};
(15,-30)="B";
(35,-30)="B'"+(4,0)*{\mathbb P^1};
"B";"B'" **\dir{-};
{\ar (10,-15)*{}; (18,-25)*{}};
{\ar (40,-15)*{}; (32,-25)*{}};
\endxy
$$\caption{}\label{o1}
\end{figure}
Since $M$ is $3$-divisible, so is the divisor
$$M_{\alpha}' := (\nu_{\alpha})_{*}(M \vert_{D_{\alpha}}) =
\sum_{i, j} a_{ij\alpha}l_{ij\alpha}'\,\, .$$

Our $D_{\alpha}'$ belongs to No.39 in the list of \cite{OS}. In
particular, the Mordell-Weil group has a torsion element of order
$3$. Thus, there are three sections $s_{0}$, $s_{1}$ and $s_{2}$
which meet $l_{00\alpha}'$, $l_{01\alpha}'$, $l_{02\alpha}'$
respectively. On the other hand, since $M_{\alpha}' \cdot
l_{ij\alpha}'$ are divisible by $3$, the set of $3$ elements
$\{a_{i0\alpha}, a_{i1\alpha}, a_{i2\alpha}\}$ (counted with
multiplicities) is either one of $\{0, 0, 0\}$ ,$\{1,1,1\}$,
$\{2,2,2\}$, $\{0, 1, 2\}$, for each $i = 0$, $1$, $2$. Suppose
that for $i = 0$ we have $\{a_{00\alpha}, a_{01\alpha},
a_{02\alpha}\} = \{0, 1, 2\}$. Then, the same holds for $i = 1$
and $2$. This is because $s_{0}\cdot M_{\alpha}'$, $s_{1}\cdot
M_{\alpha}'$, $s_{2} \cdot M_{\alpha}'$ are all divisible by $3$.
Thus both $\vert \Lambda \cap \{(i, j, \alpha)\, \vert\, i, j =
0, 1, 2\} \vert$ and $\vert \Lambda' \cap \{(i, j, \alpha)\,
\vert\, i, j = 0, 1, 2\} \vert$ are divisible by $3$ for each
$\alpha \in \{0, 1, 2\}$. This implies the result.
\end{proof}

Now we are ready to prove $6$-divisibility of the linear from $c_{2}(X)$.
It suffices to check that $D \cdot c_{2}(X) \equiv 0\, {\rm mod}\, 6$ for $D_{ijl}$,
$E_{ijk}$ and $T_{\Lambda, \Lambda'}$ in Proposition \ref{fidx}.

We have $K_{E_{ijk}}^{2} = E_{ijk}^{3} = 9$ and $c_{2}(E_{ijk}) = 3$, as
$E_{ijk} \simeq \mathbf P^2$. Thus $E_{ijk}\cdot c_{2}(X) = -6$ by Lemma
\ref{chn}. This also implies $6$-divisibility of
$T_{\Lambda, \Lambda'}\cdot c_{2}(X)$ as both $\vert \Lambda \vert$ and
$\vert \Lambda' \vert$ are divisible by $3$.

Let us compute $D_{ijl} \cdot c_{2}(X)$. As we have observed in Lemma \ref{mult},
the surface $D_{ijl}$
is the blow up at three points of a relatively minimal rational elliptic
surface. Thus, $K_{D_{ijl}}^{2} = -3$ and $c_{2}(D_{ijl}) = 15$, and therefore,
$D_{ijl}\cdot c_{2}(X) = 18$ by Lemma \ref{chn}.

This completes the proof Theorem \ref{cbc}.

\section{Connecting $X_\phi$ and $X_T$ by projective flat deformation} \label{main}

In this section we shall prove (\ref{conn}) in Theorem
\ref{mainthm}, i.e. that $X_\phi$ and $X_T$ are connected by
projective flat deformation. By Theorem \ref{key}, this follows
from:

\begin{theorem} \label{main_prop}
There are ample divisors $H_\phi$ on $X_\phi$ and $H_T$ on $X_T$  such that
$$H_\phi \cdot c_2(X_\phi) = H_T \cdot c_2(X_T)\,\, {\rm and}\,\, H_\phi^3 = H_T^3\,\, .$$
\end{theorem}
We shall prove Theorem \ref{main_prop} in the sequel. In the
proof, we freely use the notations in Notation \ref{lst}.

\subsection{Construction of a divisor $H_\phi$ on $X_\phi$}\

Recall that $E_\zeta/\langle \zeta \rangle  \simeq \mathbb P^1$.
 Let $\overline L_i=p_i^*\mathcal O_{E_\zeta/\langle \zeta \rangle  }(1)$ and ${L_i}=\pi^*\overline L_i $. Let
$$H_\phi = - \sum_{i,j,k}E_{ijk} + x L_1 +y L_2 +z L_3\,\, ,$$
where $x, y$ and $z$ are positive integers.
\begin{lemma} \label{H-phi}

\begin{enumerate}
\item For sufficiently large number $C$, $H_\phi$ is ample on $X_\phi$
when $x > C, y > C, z >  C$.
\item $H_\phi \cdot c_2(X_\phi) = 162.$
\item $H_\phi^3 = 54 xyz - 243.$
\end{enumerate}
\end{lemma}
\begin{proof} By construction, the divisor $-\sum_{i,j,k}E_{ijk}$ is
$\pi$-ample, the divisors $\overline L_i $'s are nef
on $\overline X$ and $\overline L_1 +\overline L_2 +\overline L_3 $ is
ample on $\overline X$. This implies (1).
Note that $L_i$ is represented by a smooth abelian surface
and $E_{ijk} \simeq \mathbf P^2$. Thus, by Lemma \ref{chn},
we have $L_i \cdot c_2(X_\phi) = 0$
and $E_{ijk}\cdot c_{2}(X_\phi) = -6$. This implies (2).
Note also that
$$E_{ijk}^3 = 9\,\, , \,\, L_1 \cdot L_2 \cdot L_3 = 9\,\, , \,\,
E_{ijk} \cdot L_l = L_i^2=0\,\, ,$$
and $E_{ijk} \cdot E_{lmn} =0$ unless $(i,j,k) = (l,m,n)$.
Therefore we have
\begin{align*}
H_\phi^3 =&  (- \sum_{i,j,k}E_{ijk} )^3 + 3 (- \sum_{i,j,k}E_{ijk}  )^2 (x L_1 +y L_2 +z L_3 ) \\
          &+ 3 (- \sum_{i,j,k}E_{ijk} )(x L_1 +y L_2 +z L_3 )^2  \\
          &+ (x L_1 +y L_2 +z L_3 )^3\\
         =& - \sum_{i,j,k}E_{ijk}^3 + 0 + 0+ 6 xyz L_1 \cdot L_2 \cdot L_3  \\
         =&54 xyz - 243\,\, .
\end{align*}
\end{proof}

\subsection{Construction of a divisor $H_T$ on $X_T$}\

We recall the following commutative diagram:
$$
\xymatrix{
X_\phi \ar[r]^{\pi} \ar@{.>}[d] \ar[dr]&\overline X  \ar[r]^{p_3} \ar[d]^{p_{12}}& E_\zeta/\langle \zeta \rangle  \simeq \mathbf P^1\\
X_T \ar[r]_{p_T} & B   \ar[r]_{g_1}&E_\zeta/\langle \zeta \rangle \simeq  \mathbf P^1\,\, .
}$$

Let $l_i' = g_1^{-1}(\overline Q_i)$ and $M_i = \overline
{p_{T}^{-1} ( l_i' \setminus \{\overline Q_{i0}, \overline Q_{i1},
\overline Q_{i2}\}})$ ($i = 0$, $1$, $2$). Then $M_i$ is a
relatively minimal rational elliptic surface. We denote a general
smooth fiber of the fibration $M_i \longrightarrow \mathbf P^1$ by
$f_{M_i}$. By construction, $M_{i}$ has 3 singular fibers of
Kodaira type $IV$:
$$
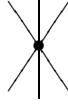
\xy
(0,0)*{\bullet};
(0,-7)*{}="A";
(0,7)*{}="A'";
(-4,6)*{}="B";
(4,-6)*{}="B'";
(-4,-6)*{}="C";
(4,6)*{}="C'";
"A";"A'" **\dir{-};
"B";"B'" **\dir{-};
"C";"C'" **\dir{-};
\endxy
$$
See Figure \ref{m-figure} for $M_i$ and the way $M_{i}$ intersects with
$F_{ijk}$'s.
\begin{figure}[h]
$$ \xy
(-40,-10)*{}="A";
(-36,-10)*{}="A'";
(-24,-10)*{}="A''";
(-20,-10)*{}="B";
(-20,10)*{}="C";
(-40,10)*{}="D"+(-3,2)*{F_{i00}};
"A";"A'" **\dir{-};
"A'";"A''" **\dir{..};
"A''";"B" **\dir{-};
"B";"C" **\dir{-};
"C";"D" **\dir{-};
"D";"A" **\dir{-};
(-36,-16)*{}="A";
(-30,-10)*{}="A'";
(-25,-5)*{}="B";
(-25, 10)*{}="B'";
(-25,15)*{}="C"+(1,3)*{F_{i01}};
(-36,4)*{}="D";
"A";"A'" **\dir{-};
"A'";"B" **\dir{..};
"B";"B'" **\dir{..};
"B'";"C" **\dir{-};
"C";"D" **\dir{-};
"D";"A" **\dir{-};
(-35,-5)*{}="A";
(-30,-10)*{}="A'";
(-24,-16)*{}="B";
(-24,4)*{}="C";
(-35,15)*{}="D"+(-1,3)*{F_{i02}};
(-35,10)*{}="D'";
"A";"A'" **\dir{..};
"A'";"B" **\dir{-};
"B";"C" **\dir{-};
"C";"D" **\dir{-};
"D";"D'" **\dir{-};
"D'";"A" **\dir{..};
(-30,-10)*{}; (-30,10)*{} **\dir{-};
{\ar@/^.5pc/(-53,-20)+(-7, 3)*{\stackrel{\text{\large {ruling in}}} { F_{i01} \simeq \mathbb F_1}  };(-33,-12)*{}};
{\ar@/^-.5pc/(-23,-26)+(0,-7)*{ \text{(-2) curve in } M_i};(-31,-12)*{}};
(-10,-10)*{}="A";
(-6,-10)*{}="A'";
(6,-10)*{}="A''";
(10,-10)*{}="B";
(10,10)*{}="C";
(-10,10)*{}="D"+(-3,2)*{F_{i10}};
"A";"A'" **\dir{-};
"A'";"A''" **\dir{..};
"A''";"B" **\dir{-};
"B";"C" **\dir{-};
"C";"D" **\dir{-};
"D";"A" **\dir{-};
(-6,-16)*{}="A";
(0,-10)*{}="A'";
(5,-5)*{}="B";
(5, 10)*{}="B'";
(5,15)*{}="C"+(1,3)*{F_{i11}};
(-6,4)*{}="D";
"A";"A'" **\dir{-};
"A'";"B" **\dir{..};
"B";"B'" **\dir{..};
"B'";"C" **\dir{-};
"C";"D" **\dir{-};
"D";"A" **\dir{-};
(-5,-5)*{}="A";
(0,-10)*{}="A'";
(6,-16)*{}="B";
(6,4)*{}="C";
(-5,15)*{}="D"+(-1,3)*{F_{i12}};
(-5,10)*{}="D'";
"A";"A'" **\dir{..};
"A'";"B" **\dir{-};
"B";"C" **\dir{-};
"C";"D" **\dir{-};
"D";"D'" **\dir{-};
"D'";"A" **\dir{..};
(0,-10)*{}; (0,10)*{} **\dir{-};
(20,-10)*{}="A";
(24,-10)*{}="A'";
(36,-10)*{}="A''";
(40,-10)*{}="B";
(40,10)*{}="C";
(20,10)*{}="D"+(-3,2)*{F_{i20}};
"A";"A'" **\dir{-};
"A'";"A''" **\dir{..};
"A''";"B" **\dir{-};
"B";"C" **\dir{-};
"C";"D" **\dir{-};
"D";"A" **\dir{-};
(24,-16)*{}="A";
(30,-10)*{}="A'";
(35,-5)*{}="B";
(35, 10)*{}="B'";
(35,15)*{}="C"+(1,3)*{F_{i21}};
(24,4)*{}="D";
"A";"A'" **\dir{-};
"A'";"B" **\dir{..};
"B";"B'" **\dir{..};
"B'";"C" **\dir{-};
"C";"D" **\dir{-};
"D";"A" **\dir{-};
(25,-5)*{}="A";
(30,-10)*{}="A'";
(36,-16)*{}="B";
(36,4)*{}="C";
(25,15)*{}="D"+(-1,3)*{F_{i22}};
(25,10)*{}="D'";
"A";"A'" **\dir{..};
"A'";"B" **\dir{-};
"B";"C" **\dir{-};
"C";"D" **\dir{-};
"D";"D'" **\dir{-};
"D'";"A" **\dir{..};
(30,-10)*{}; (30,10)*{} **\dir{-};
(-50,-25)*{}="A";
(55,-25)*{}="B"+(3,10)*{M_i};
(45,0)*{}="C";
(41,0)*{}="C'";
(19,0)*{}="C''";
(11,0)*{}="C'''";
(-11,0)*{}="C''''";
(-19,0)*{}="C'''''";
(-41,0)*{}="C''''''";
(-45,0)*{}="D";
"A";"B" **\dir{-};
"B";"C" **\dir{-};
"C";"C'" **\dir{-};
"C''";"C'''" **\dir{-};
"C''''";"C'''''" **\dir{-};
"C''''''";"D" **\dir{-};
"D";"A" **\dir{-};
(5, -30)*{}="A"+(2, -3)*{f_{M_i}};
(6, -30)*{}="B";
(-17, -25);(-15, 0)**\crv{(-18, -10)&(-14,-20)};
(18, -25);(15, 0)**\crv{(14, -10)&(18,-20)};
{\ar@/^-.5pc/"A";(-14.5,-12)*{}};
{\ar@/^.5pc/"B";(15.5,-12)*{}};
\endxy $$
\caption{}\label{m-figure}
\end{figure}

Let $\overline S_j = (p_3^{*}(\overline Q_j))_{\rm red}$ on
$\overline X $ and $S_j$ be the proper transformation of
$\overline S_j$ on $X_T$ ($j = 0$, $1$, $2$). Then $S_j$ is a
(non-relatively minimal) rational elliptic surface with three
singular fibers (denote them by $\eta_1, \eta_2$ and $\eta_3$)
that are composed of one $(-1)$-curve of multiplicity $3$ and
three $(-3)$-curves; $\eta_i=  \alpha_i + \beta_i + \gamma_i + 3
\delta_{i}$. See Figure \ref{Sf}. We denote by $f_{S_j}$ a general
smooth fiber of the fibration $S_j \longrightarrow \mathbf P^1$.
\begin{figure}[h]
$$ \xy
(10,0)*{\bullet}+(-2,-8)*{\alpha_i}+(-7,8)*{3 \delta_{i}};
(20,0)*{\bullet}+(-2,-8)*{\beta_i};
(30,0)*{\bullet}+(-2,-8)*{\gamma_i}; (5,0)*{}; (35,0)*{}
**\dir{-}; (5,.3)*{}; (35,.3)*{} **\dir{-}; (5,-.3)*{};
(35,-.3)*{} **\dir{-}; (8,-5)*{}; (12,5)*{} **\dir{-}; (18,-5)*{};
(22, 5)*{} **\dir{-}; (28,-5)*{}; (32,5)*{} **\dir{-};
(10,0)*{}="A"+(2,7)*{(-3)}; (20,0)*{}="B"+(2,7)*{(-3)};
(30,0)*{}="C"+(2,7)*{(-3)}; (40,0)*{}="D"+(1,0)*{(-1)};
\endxy
$$
\caption{}\label{Sf}
\end{figure}

See Figure \ref{SM} for the configuration of $S_j$, $M_i$ and
$F_{\alpha \beta \gamma}$'s.

\begin{figure}[h]
$$ \xy
(-40,-10)*{}="A";
(-37.5,-10)*{}="A'";
(-20,-10)*{}="B";
(-20,.5)*{}="B'";
(-20,10)*{}="C";
(-40,10)*{}="D"+(-3,2)*{F_{i00}};
"A";"A'" **\dir{-};
"A'";"B" **\dir{..};
"B";"B'" **\dir{..};
"B'";"C" **\dir{-};
"C";"D" **\dir{-};
"D";"A" **\dir{-};
(-36,-16)*{}="A";
(-30,-10)*{}="A'";
(-25,-5)*{}="B";
(-25, 10)*{}="B'";
(-25,15)*{}="C"+(1,3)*{F_{i01}};
(-36,4)*{}="D";
"A";"A'" **\dir{..};
"A'";"B" **\dir{..};
"B";"B'" **\dir{..};
"B'";"C" **\dir{-};
"C";"D" **\dir{-};
"D";"A" **\dir{-};
(-35,-5)*{}="A";
(-24,-16)*{}="B";
(-24,0.3)*{}="B'";
(-24,4)*{}="C";
(-35,15)*{}="D"+(-1,3)*{F_{i02}};
(-35,10)*{}="D'";
"A";"B" **\dir{..};
"B";"B'" **\dir{..};
"B'";"C" **\dir{-};
"C";"D" **\dir{-};
"D";"D'" **\dir{-};
"D'";"A" **\dir{..};
(-30,-10)*{}="E";
(-30,1.3)*{}="E'";
(-30,10)*{}="F";
"E";"E'" **\dir{..};
"E'";"F" **\dir{-};
(-10,-10)*{}="A";
(10,-10)*{}="B";
(10,0.3)*{}="B'";
(10,10)*{}="C";
(-10,10)*{}="D"+(-3,2)*{F_{i10}};
(-10,3)*{}="D'";
"A";"B" **\dir{..};
"B";"B'" **\dir{..};
"B'";"C" **\dir{-};
"C";"D" **\dir{-};
"D";"D'" **\dir{-};
"D'";"A" **\dir{..};
(-6,-16)*{}="A";
(5,-5)*{}="B";
(5, 10)*{}="B'";
(5,15)*{}="C"+(1,3)*{F_{i11}};
(-6,4)*{}="D";
"A";"B" **\dir{..};
"B";"B'" **\dir{..};
"B'";"C" **\dir{-};
"C";"D" **\dir{-};
"D";"A" **\dir{-};
(-5,-5)*{}="A";
(6,-16)*{}="B";
(6,.6)*{}="B'";
(6,4)*{}="C";
(-5,15)*{}="D"+(-1,3)*{F_{i12}};
(-5,10)*{}="D'";
"A";"B" **\dir{..};
"B";"B'" **\dir{..};
"B'";"C" **\dir{-};
"C";"D" **\dir{-};
"D";"D'" **\dir{-};
"D'";"A" **\dir{..};
(0,-10)*{}="E";
(0,2.3)*{}="E'";
(0,10)*{}="F";
"E";"E'" **\dir{..};
"E'";"F" **\dir{-};
(20,-10)*{}="A";
(36,-10)*{}="A''";
(40,-10)*{}="B";
(40,10)*{}="C";
(20,10)*{}="D"+(-3,2)*{F_{i20}};
(20,2)*{}="D'";
"A";"A''" **\dir{..};
"A''";"B" **\dir{-};
"B";"C" **\dir{-};
"C";"D" **\dir{-};
"D";"D'" **\dir{-};
"D'";"A" **\dir{..};
(24,-16)*{}="A";
(30,-10)*{}="A'";
(26,-14.3)*{}="A''";
(35,-5)*{}="B";
(35, 10)*{}="B'";
(35,15)*{}="C"+(1,3)*{F_{i21}};
(24,4)*{}="D";
"A'";"B" **\dir{..};
"A";"A''" **\dir{..};
"A''";"A'" **\dir{-};
"B";"B'" **\dir{..};
"B'";"C" **\dir{-};
"C";"D" **\dir{-};
"D";"A" **\dir{-};
(25,-5)*{}="A";
(30,-10)*{}="A'";
(36,-16)*{}="B";
(36,4)*{}="C";
(25,15)*{}="D"+(-1,3)*{F_{i22}};
(25,10)*{}="D'";
"A";"A'" **\dir{..};
"A'";"B" **\dir{-};
"B";"C" **\dir{-};
"C";"D" **\dir{-};
"D";"D'" **\dir{-};
"D'";"A" **\dir{..};
(30,-10)*{}; (30,10)*{} **\dir{-};
(-37.5,-15.3)*{}="A";
(26,-15)*{}="B";
(26,4.9)*{}="C";
(-37,4.7)*{}="D";
(-37.5,-35.3)*{}="A'";
(26,-35)*{}="B'";
(-37.5, -25)="A''";
(26, -25)="B''";
(-37.5,-30.3)*{}="E";
(26,-30)*{}="F";
"A";"B"  **\crv{(-20, -25) & (-6,-7) & (7, -25)};
"B";"C" **\dir{-};
"D";"C" **\crv{(-20, -5) & (-6, 13) & (7, -5)};
"D";"A" **\dir{-};
"A'";"B'"  **\crv{(-20, -45) & (-6,-27) & (7, -45)};
"E";"F"  **\crv{(-20, -40) & (-6,-28) & (7, -35)};
"A";"A''" **\dir{..};
"A''";"A'" **\dir{-};
"B''";"B'" **\dir{-};
"B";"B''" **\dir{..};
(-50,-25)*{}="A";
(55,-25)*{}="B"+(3,10)*{M_i};
(45,0)*{}="C";
(41,0)*{}="C'";
(19,0)*{}="C''";
(11,0)*{}="C'''";
(-11,0)*{}="C''''";
(-19,0)*{}="C'''''";
(-41,0)*{}="C''''''";
(-45,0)*{}="D";
"A";"B" **\dir{-};
"B";"C" **\dir{-};
"C";"C'" **\dir{-};
"C''''''";"D" **\dir{-};
"D";"A" **\dir{-};
(-47,-0) = "L"+(-15,0)*{\stackrel{ \text{ \large (-1) sing. fiber} } {  \text{in } S_1 \text{ and section}}}+(-2, -6)*{{\text{ of } M_i  \text{ in }M_i}};
(30, -30)="S"+(3,0)*{S_1};
(-45, -40)="F"+(-3,0)*{f_{S_1}};
{\ar@/^.5pc/"S";(15,-10)*{}};
{\ar@/^1pc/"L";(-30,-18)};
{\ar@/^-1pc/"F";(-30,-34.7)};
(-10,-45) = "K"+(0,-5)*{\stackrel{ \text{ \large 1-section with}\,\,\,} { s^2 =1  \text{ in } F_{i11} \text{ and}}}+(-2, -6)*{{\text {(-3) curve in }}S_1};
{\ar@/^1pc/"K";(-6.3,-13)};
\endxy $$
\caption{}\label{SM}
\end{figure}

\begin{lemma}\label{l_H_T}
The following divisor is $p_{T}$-ample;
$$3(M_0 + M_1 + M_2) + S_0 + S_1 + S_2\,\, .$$
\end{lemma}
\begin{proof}
Since $S_j$ are sections of $p_{T}$ over $B \setminus \{ \overline
Q_{\alpha \beta}\, {\big|}\, \alpha, \beta =0,1,2 \}$, we only
need to check that $3(M_0 + M_1 + M_2) + S_0 + S_1 + S_2$ is ample
on $F_{\alpha \beta \gamma}$'s. This, however, follows from the
fact that:
$$\left(3(M_0 + M_1 + M_2) + S_0 + S_1 + S_2)\right)|_{F_{\alpha \beta \gamma}}= 3f + s\,\, ,$$
where $f$ is the ruling of the ruled surface $F_{\alpha \beta \gamma}
\simeq \mathbf F_{1}$
and $s$ is a positive section (with $s^{2} = 1$). (See Figure \ref{m-figure},  \ref{SM}).
\end{proof}
Let $\overline A_k$ be a general fiber of $g_{k} : B
\longrightarrow E_\zeta/\langle \zeta \rangle$ and let $A_k =
p_{T}^{*}\overline A_k$ ($k = 1$, $2$). $\overline A_k$ is an
elliptic curve and $A_k$ is an abelian surface. Note that
$\overline A_k \in \vert g_k^* \mathcal O_{E_\zeta/\langle \zeta
\rangle}(1) \vert$. Put
$$H_T:=3(M_0 + M_1 + M_2) + S_0 + S_1 + S_2 + a A_1 +b A_2\,\, .$$
Here $a$ and $b$ are positive integers.
\begin{lemma}\label{H-T}

\begin{enumerate}
\item For sufficiently large number $C$, $H_T$ is ample on $X_T$ when $a > C$,
$b > C$.
\item $H_T \cdot c_2(X_T) = 162$.
\item $H_T^3 = 18ab - 27 b -333$.
\end{enumerate}
\end{lemma}
\begin{proof}
Since $\overline A_1, \overline A_2$ are nef and $\overline A_1 + \overline A_2$
is ample on $B$, the first assertion follows from Lemma
\ref{l_H_T}. By using Lemma\ref{chn}, we compute that
\begin{align*}
M_i \cdot c_2(X_T) =& -K_{M_i}^2 + c_2(M_i) = 12\\
S_j \cdot c_2(X_T) =& -K_{S_j}^2 + c_2(S_j) = 18\\
A_k \cdot c_2(X_T) =& -K_{A_k}^2 + c_2(A_k) =  0\,\, .
\end{align*}
This implies the second assertion.
For the third one, we first expand $H_T^3$ as:
\begin{align*}
H_T^3 =& \left( 3(M_0 + M_1 + M_2) + S_0 + S_1 + S_2 \right)^3 &(= Q_1)\\
       &+\left( 3(M_0 + M_1 + M_2) + S_0 + S_1 + S_2 \right)^2 (aA_1 + b A_2) &(= Q_2)\\
       &+ \left( 3(M_0 + M_1 + M_2) + S_0 + S_1 + S_2 \right) (aA_1 + b A_2)^2 &(= Q_3)\\
       &+ (aA_1 + b A_2)^3 &(=Q_4)\\
       =&Q_1 + Q_2 + Q_3 + Q_4.&
\end{align*}
We compute $Q_{1}$, $Q_{2}$, $Q_{3}$, $Q_{4}$ separately.
\begin{description}
\item[$Q_1$] Note that
\begin{align*}
&S_j^3 = K_{S_j}^2 = -3\\
&M_i^3 = K_{M_i}^2 = 0\\
&S_i \cdot S_j = M_i \cdot M_j =0 \text{ for } i \neq j\\
&M_i^2 \cdot S_j = (M_i|_{S_j})^2 = -1 \\
&M_i \cdot S_j^2 = (S_j|_{M_i})^2 = -1\,\, .
\end{align*}
With these, we have $Q_1 = -333$. \item[$Q_2$] We observe that
\cen{$-M_i|_{M_i} \sim -K_{M_i} \sim A_1|_{M_i} \sim f_{M_i}$ and
$f_{M_i} \cdot A_2= 3$.} From this we have $M_i^2 \cdot (aA_1 + b
A_2) =-3b$. Note also that \cen{$A_1|_{S_j} \sim f_{S_j}$.} It
follows that
\begin{align*}
M_i \cdot S_j \cdot (a A_1 + b A_2) =& M_i|_{S_j} \cdot (a A_1 + b A_2)|_{S_j}\\
                                    =&b(\delta_{i}\cdot A_2)\\
                                    =&b\,\, .
\end{align*}
See also Figure \ref{Sf}. Finally,
\begin{align*}
S_j^2 \cdot (a A_1 + b A_2) =& K_{S_j} \cdot (a A_1 + b A_2)|_{S_j}\\
                            =& b(-f_{S_j} + \delta_{1} + \delta_2 + \delta_3) \cdot A_2\\
                            =& b(-3 + 1 + 1 +1)\\
                            =&0.
\end{align*}
Thus we have $Q_2 = -27b$.
\item[$Q_3$] Note that $(aA_1 + bA_2)^2 = 6ab (\text{fiber of }p_{T}) $. So we have
\begin{align*}
&M_i\cdot(a A_1 + b A_2)^2 = 0\\
&S_j \cdot(a A_1 +b A_2)^2 = 6ab
\end{align*}
Thus, $Q_3 = 18ab$.
\item[$Q_4$] Clearly, $Q_4=0$.
\end{description}
With all these, we obtain $H_T^3 = Q_1 + Q_2 + Q_3 + Q_4 = 18ab -27b -333$.
\end{proof}

\subsection{Synthesis}\

Now we are ready to prove Theorem \ref{main_prop}. By Lemma
\ref{H-phi} and Lemma \ref{H-T}, the divisors $H_\phi$ and $H_T$
are ample on $X_\phi$ and $X_T$ respectively when $x, y, z$ and
$a, b$ are greater than some sufficiently large $C$. So, it
suffices to find integers $x, y, z$ and $a, b$ greater than any
given positive integer $C$ that satisfy the following equations:
$$162=H_\phi \cdot c_2(X_\phi) = H_T \cdot c_2(X_T)=162$$
and
$$54 xyz - 243=H_\phi^3 = H_T^3 =18ab - 27 b -333.$$
The first one poses no condition on $x$, $y$, $z$, $a$, $b$, and the
second one is simplified to:
\begin{align}\label{key-eqn}
6 xyz  =2ab - 3 b -10.
\end{align}
For a given positive integer $C$, let
$$x = 12 C^2 -6\,\, ,\,\, y = z = 2 C\,\, ,\,\, a = 6 C^2 +1\,\, ,\,\,
b = 24C^2-10\,\, .$$ Then $x,y,z$ and $a, b$ are integers which
are greater than $C$ and satisfy the above equation
(\ref{key-eqn}).

This completes the proof of Theorem \ref{main_prop}.


\end{document}